\documentclass[envcountsame]{mmnp}
\usepackage{amsmath,amssymb}
\usepackage{amsfonts}
\usepackage[english]{babel}

\usepackage[latin1]{inputenc}
\usepackage{setspace}
\usepackage{dsfont}
\usepackage{graphics,subfigure,xcolor}

\usepackage{tikz,pgflibraryshapes,xkeyval}
\usepackage{caption}

\newcommand {\col}{{\mathrm{col}}}

\newfont{\pseudocode}{cmtt10}

\newcommand{\pd}{{\partial}}

\newcommand {\ess}{{\mathrm{ess}}}

\newcommand{\Real}{\mathbb{R}}

\newcommand{\rank}{\mathrm{rank}}

\renewcommand{\col}[1]{\mathrm{col}\left(#1\right)}

\newcommand{\beq}{\begin{equation}}
\newcommand{\eeq}{\end{equation}}

\newcommand{\bbm}{\begin{bmatrix}}
\newcommand{\ebm}{\end{bmatrix}}

\newcommand{\bpm}{\begin{pmatrix}}
\newcommand{\epm}{\end{pmatrix}}

\newcommand{\bit}{\begin{itemize}}
\newcommand{\eit}{\end{itemize}}

\newcommand{\ben}{\begin{enumerate}}
\newcommand{\een}{\end{enumerate}}

\newcommand{\barr}{\begin{array}}
\newcommand{\earr}{\end{array}}

\idline{Vol. , No. }{1}
\doi{10.1051/mmnp/2016xxxx}

\numberwithin{equation}{section}

\newcounter{assumption}[section]
\renewcommand{\theassumption}{\thesection.\arabic{assumption}}
\newenvironment{assumption}[1][]{\refstepcounter{assumption}\par\medskip
   \noindent \textbf{Assumption~\theassumption.} {\textit{#1}} \rmfamily}{\medskip}

\begin{document}

\title{Fast Sampling of Evolving Systems with Periodic Trajectories}

\author{I.Yu. Tyukin\inst{1,2}\thanks{\email {I.Tyukin@le.ac.uk}}, A.N. Gorban\inst{1}, T.A. Tyukina\inst{1}, J. Mohammed Al-Ameri\inst{1}, Yu.A. Korablev\inst{2}}

\vspace{0.5cm}

\institute{\inst{1} {University of Leicester, Department of Mathematics, United Kingdom}\\
\inst{2} {Saint-Petersburg State Electrotechnical University,\\ Department of Automation and Control Processes, Russia}}



\titlerunning{Fast Sampling of Evolving Systems}

\abstract{We propose a novel method for fast and scalable evaluation of periodic solutions of systems of ordinary differential equations for a given set of parameter values and initial conditions. The equations governing the system dynamics are supposed to be of a special class, albeit admitting nonlinear parametrization and nonlinearities.  The method enables to represent a given periodic solution as sums of computable integrals and functions that are explicitly dependent on parameters of interest and initial conditions. This allows invoking parallel computational streams in order to increase speed of calculations. Performance and practical implications of the method  are illustrated with examples including classical predator-prey system and models of neuronal cells.
}

\maketitle

\section*{Notation}

\begin{itemize}
\item[$\bullet$] Symbol $\|\cdot\|$ stands for the Euclidian norm.
\item[$\bullet$] Symbol $\Real$ denotes the field of real numbers, and $\Real_{\geq 0}=\{x\in\Real \ | \ x\geq 0\}$.
\item[$\bullet$] Symbols $I_n$, $0_n$ are the identity and zero $n\times n$ matrices, respectively.
\item[$\bullet$] By $\mathcal{K}$ we denote   the set of all strictly increasing
continuous functions $\kappa: \Real_{\geq 0}\rightarrow
\Real_{\geq 0}$ such that $\kappa(0)=0$.
\item[$\bullet$] Consider a non-autonomous system $\dot{x}=f(x,p,t,u(t))$,  where $f:\Real^n\times\Real^d\times\Real\times\Real^l\rightarrow\Real^n$, $u:\Real\rightarrow\Real^l$
are continuous, $p\in\Real^d$ is the vector of parameters, and $f(\cdot,p,t,u)$ is locally Lipschitz; $x( \cdot \ ;t_0,x_0,p,[u])$ stands for
the unique maximal solution of the initial value problem: $x(t_0;t_0,x_0,p,[u])=x_0$. In cases when no confusion arises, we will refer to these solutions as $x(\cdot;t_0,x_0,[u])$,  $x(\cdot;x_0,[u])$,  or simply $x(\cdot)$. Solutions of the initial value problem above at $t$ are denoted as $x(t;t_0,x_0,p,[u])$,
$x(t;t_0,x_0,[u])$,  $x(t;x_0,[u])$,  or $x(t)$ respectively.
\item[$\bullet$] Let $f:\Real\rightarrow\Real^n$, then $\|f(\tau)\|_{\infty,[t_0,t_0+T]}$ denotes the uniform norm of $f(\cdot)$ on $[t_0,t_0+T]$:  $\|f(\tau)\|_{\infty,[t_0,t_0+T]}=\ess \sup\{\|f(t)\|,t \in [t_0,t_0+T]\}$.
\end{itemize}

\setcounter{equation}{0}

\section{Introduction}

Evolutionary development and change is a wide-spread phenomenon that is inherent to biological systems. It can be viewed as a change of hereditary characteristics of an organism over successive generations. In mathematical biology hereditary characteristics may be imagined, in loose terms, as a model variable determining dynamics of populations that changes over time. Alternations of individual traits often occur at a microscopic level; the effects, however, can be felt at the level of single organisms, their groups as well as their populations. Sometimes these effects will result in a gradual change of fitness in a near future, but on a longer time scale and near crises evolutionary changes may lead to significant reduction of population and even extinction. The question is therefore if it is possible to infer minute evolutionary changes fast enough to predict potential catastrophes in future?

In mathematical modelling terms evolutionary alternations could be reflected in parameter changes of equations determining dynamics of species, organisms or populations. A simple example could be the changes of predation or reproduction rates in predator-prey models with saturation \cite{JCAM:Chen}, \cite{CSB:Chen} (cf \cite{Ecology:Alonso})
\begin{equation}\label{eq:predator_prey}
\begin{split}
\dot{x}&= p_1 x \left(1 - \frac{x}{p_2} \right) - \frac{p_3 z x}{p_4  +  x}\\
\dot{z}& =  \frac{p_5 z x}{p_4  + x} - p_6 z,
\end{split}
\end{equation}
where $x$, $z$ are population densities of prey and predator, respectively, and $p=(p_1,\dots,p_6)$ are parameters that are subjected to evolutionary modifications. Models of this type are known to exhibit broad range of qualitatively different dynamics, depending on parameters \cite{CSB:Chen}, including stable limit cycles and relaxations to equilibria. Hence the question of monitoring and predicting potential consequences of minor evolutionary changes in relevant systems can be posed as a problem of fast parameter and state estimation of (\ref{eq:predator_prey}) from available measurement data. The data, in turn, may be restricted to few samples of either $x(t)$ or $z(t)$.

In a more generalized framework we will assume that processes of interest can be described by  systems of nonlinear ordinary differential equations
\begin{equation}\label{eq:system0}
\dot{x}=f(x,p,t), \ x(t_0)=x_0,
\end{equation}
where $f:\Real^n\times\Real^k\times\Real\rightarrow\Real^n$ is \ continuous and locally Lipschitz wrt the variable $x$ function, and $p$ is the vector of unknown parameters. Let $[t_0,t_0+T]$ be an interval on which the solution $x(\cdot;t_0,x_0,p)$ of (\ref{eq:system0}) is defined. Let us further suppose that the system's state, $x(t;t_0,x_0,p)$, is not accessible for direct observation at any $t\in[t_0,t_0+T]$. We will, however, assume that the values
\[
h(t,x(t;t_0,x_0,p)), \ h:\Real\times\Real^n\rightarrow\Real
\]
are available for every $t\in[t_0,t_0+T]$. 
The problem is now how to find  $p'\in \Real^k$, $x_0'\in\Real^n$ such that
\begin{equation}\label{eq:measurement:match}
\begin{array}{c}
h(t,x(t;t_0,x_0,p))=h(t,x(t;t_0,x_0',p')) \ \mbox{for all}  \  t\in [t_0,t_0+T].
\end{array}
\end{equation}

This is a standard inverse problem, and many methods for finding solutions to this  problem have
been developed to date (sensitivity functions \cite{SIAM:Ident}, splines \cite{Brewer:2008}, interval analysis \cite{Automatica:2008:Johnson},  adaptive observers \cite{Marino92},\cite{Besancon:2000}, \cite{Automatica:Farza:2009}, \cite{Automatica:Grip:2010},\cite{Tyukin:2011},\cite{Tyukin_2012:arx}, \cite{MMNP:2010:tyukin}  and particle filters and Bayesian inference methods \cite{Abarbanel:2009}). Despite these methods are based on different mathematical frameworks, they share a common feature:
one is generally required to repeatedly  find numerical solutions of nonlinear ordinary differential equations (ODEs) over given intervals of time (solve the direct problem).

In particular, if the right-hand side of (\ref{eq:system0}) is differentiable with respect to $x$ and $p$ then in the framework of sensitivity functions at each query point $\hat{p},\hat{x}_0$ one needs to solve the following initial value problem:
\begin{equation}\label{eq:sensitivity}
\begin{split}
\dot{x}&=f(t,x,\hat{p})\\
\dot{\delta}_p&=\frac{\pd f}{\pd x}[t,x(t;t_0,\hat{x}_0,\hat{p}),\hat p] \delta_p + \frac{\pd f}{\pd p}[t,x(t;t_0,\hat{x}_0,\hat{p}),\hat p]\\
\dot{\delta}_0&= \frac{\pd f}{\pd x} [t,x(t;t_0,\hat{x}_0,\hat{p}),\hat p] \delta_0\\
x(t_0)&=\hat{x}_0, \ \delta_p(t_0)=0_n, \ \delta_0(t_0)=I_n,
\end{split}
\end{equation}
where $\delta_p(t)$ is the $n\times k$ matrix of partial derivatives of solution $x(t;t_0,\hat{x}_0,\hat{p})$ with respect to $\hat{p}$, and $\delta_0(t)$ is the $n\times n$ matrix of partial derivatives of $x(t;t_0,\hat{x}_0,\hat{p})$ with respect to $\hat{x}_0$.

Notwithstanding the plausibility of numerical integration of (\ref{eq:sensitivity}) or (\ref{eq:system0}) in algorithms for state and parameter estimation, this operation is an inherently sequential process. The amount of time needed for direct numerical integration of these systems over a grid of $N$ points is at least $O(N)$. This constrains computational scalability of the problem, and as a result it imposes limitations on the time required to derive a solution.

In order to overcome this limitation we propose  to cast the inverse problem above in an alternative, integral form in which the model output is defined as a combination of indefinite integrals with known, explicit and computable kernels, possibly dependent on $x_0,p$, and explicit functions of $x_0,p$. The advantage of such integral formulations is that their computations are scalable and can be performed using conventional prefix sum algorithms of which the execution time is of order $O(\log_a(N))$, $a=2,3\dots$. The latter option compares favorably with $O(N)$.

Hence, instead of finding numerical solutions of the initial value problem (\ref{eq:system0}) and matching the results to observed data, e.g. as (\ref{eq:measurement:match}), we will search for a representation of the problem as
\begin{equation}\label{eq:integral_formulation_general}
\begin{split}
&y(t)-F(p,x_0,t)=0, \mbox{for all} \ t\in[t_0,t_0+T], \\
&F(p,x_0,t)=R(p,x_0)+\int_{t_0}^t g(t,\tau,p,x_0)d\tau,
\end{split}
\end{equation}
where $R:\Real^k\times\Real^n\rightarrow\Real^n$, $g:\Real\times\Real\times\Real^k\times\Real^n\rightarrow \Real^d$ and $y:\Real\rightarrow \Real^d$ are functions that are explicitly computable from measurement data. We additionally require that if $p',x_0'$ is a solution of (\ref{eq:integral_formulation_general}) then it is also a solution of (\ref{eq:measurement:match}) and vise versa. Finally, if the values of $h(t,x(t;t_0,x_0,p))$ are available only at a finite set of points $\{t_i\}$, $t_i\in[t_0,t_0+T]$ we will assume that they are informative enough so that the functions $R$, $g$ can be reconstructed or approximated with acceptable accuracy, and corresponding representations $F(p,x_0,t_i)$ in (\ref{eq:integral_formulation_general}) can still be  computed with sufficient precision.

In the next sections we  specify a class of systems for which such representation is possible. This class of systems is not as general as (\ref{eq:system0}) but is relevant enough in modelling applications. In Section \ref{sec:problem_statement} we define this class of systems and present general technical assumptions. This is followed by presentation of main results in Section \ref{sec:main_results}. The results are based on the periodicity assumption we impose on the data and also on known facts from the theory of adaptive observes \cite{Marino92}, \cite{Lorea_2002}. In Section \ref{sec:example} we illustrate the approach with two examples as well as provide an assessment of scalability of the approach in the problems considered.

\setcounter{equation}{0}

\section{Problem Formulation}\label{sec:problem_statement}

Consider the following class of systems
\begin{equation}\label{eq:1}
\begin{split}
\dot{x}&=A(\theta) x  + \Psi(y,t)\theta + v(y,q,\lambda,t)\\
\dot{q}&= P(y,\lambda,t) q + w(y,\lambda,t)\\
y&=C^{T}x, \ x(t_0)=x_0, \ q(t_0)=q_0,
\end{split}
\end{equation}
where $(x,q)$,  $x\in\Real^n$, $q\in\Real^d$ is the state vector,  $\theta\in\Real^m$, $\lambda\in\Real^p$ are parameters, $A(\theta)$ is an $n\times n$ real matrix, possibly dependent on $\theta$, and $C\in\Real^n$, $C=\col{1,0,\dots,0}$.  We assume that the following hold for (\ref{eq:1}):

\begin{assumption}[(General assumptions on (\ref{eq:1}))]
 $\ \ \ \ \ \ \ \ \ \ $
\begin{itemize}
\item[A1)] the solution of (\ref{eq:1}) is defined on the interval $[t_0,t_0+T]$ (for some $T>0$, possibly dependent on $t_0$);
\item[A2)] the pair $A(\theta),C^T$, is  observable, that is
\[
\rank\left(\begin{array}{c}
                    C^{T}\\
                    C^T A(\theta)\\
                     \vdots\\
                     C^{T} A^{n-1}(\theta)
            \end{array}\right)=n;
\]
\item[A3)] $P(y,\lambda,t)$ and $\Psi(y,t)$ are  $d\times d$ and $n\times m$ real matrices of which the entries are continuous and differentiable functions; $P(y,\lambda,t)$ is diagonal:
    \[
    P(y,\lambda,t)=\mathrm{diag}\left(\alpha_1(y,\lambda,t),\dots,\alpha_d(y,\lambda,t) \right);
    \]
\item[A4)] $v:\Real\times\Real^d\times\Real^p\times\Real\rightarrow \Real^n$, $w:\Real\times\Real^p\times\Real\rightarrow\Real^d$ are continuous and differentiable functions.
\item[A5)] Exact values of parameters  $\theta$, $\lambda$ are unknown.
\end{itemize}
\end{assumption}

Since the pair $A(\theta),C^T$ is observable there always is a coordinate transform, possibly  dependent on $t,\theta$, \cite{Marino92} rendering (\ref{eq:1}) into the following form
\begin{equation}\label{eq:system}
\begin{split}
\dot{x}&=A_0 x  + b \varphi(y,t)^{T} \tilde\theta(\theta) + \tilde{v}(y,q,\bar\lambda(\lambda,\theta),t)\\
\dot{q}&= \tilde P(y,\bar \lambda,t) q + \tilde w(y,\bar\lambda,t)\\
y&=C^{T}x, \ x(t_0)=x_0, \ q(t_0)=q_0,
\end{split}
\end{equation}
$b=(1,b_1,\dots,b_{n-1})^{T}$ is such that the polynomial $s^{n-1} + b_1s^{n-2}+\cdots+b_{n-1}$ is
Hurwitz, and $A_0=\left(\begin{array}{cc}0 & I_{n-1}\\
                           0 & 0\end{array}\right)$.
Functions
$\varphi:\Real\times\Real\rightarrow\Real^r$,
$\tilde{v}:\Real\times\Real^k\times\Real\rightarrow\Real^n$, $\tilde w: \Real\times\Real^k\times\Real\rightarrow \Real^d$, $\tilde P:\Real\times\Real^k\times\Real\rightarrow \Real^{d\times d}$ are
continuous and differentiable, $\tilde P$ is diagonal,  and $\tilde\theta\in\Real^r$,
$\bar\lambda\in\Real^k$ are parameters.

Furthermore, noticing that the variable $y$ is defined and known for all $t\in[t_0,t_0+T]$, and  $y(\cdot)$ is continuous one can express the solution $q(t;q_0,\bar\lambda,[y])$ on $[t_0,t_0+T]$ in the closed form as follows:
\[
\begin{split}
q(t;q_0,\bar\lambda,[y])&=e^{\int_{t_0}^t \tilde P(y(\tau),\bar\lambda,\tau)d\tau}q_0 +\int_{t_0}^t e^{\int_{\tau}^t \tilde P(y(s),\bar\lambda,s)ds} \tilde w(y(\tau),\bar\lambda,\tau) d\tau.
\end{split}
\]
Denoting  $\tilde{\lambda}=\mathrm{col}(\bar\lambda,q_0)$, $g(y,\tilde{\lambda},t)=\tilde{v}(y,q(t;q_0,\bar\lambda,[y]),\bar\lambda,t)$ we therefore arrive at the transformed equations of (\ref{eq:system}):
\begin{equation}\label{eq:system1}
\begin{split}
\dot{x}&=A_0 x  + b \varphi(y,t)^{T} \tilde\theta(\theta) + g(y,\tilde{\lambda},t)\\
y&=C^{T}x, \ x(t_0)=x_0.
\end{split}
\end{equation}

The core problem we are interested in (\ref{eq:system1}) is as follows:
\begin{prblm}\label{prob:inferrence} Let (\ref{eq:system1}) be given, and its solutions are defined on $[t_0,t_0+T]$. Suppose that all functions in the right-hand side of (\ref{eq:system}) are known, but  true values of $x_0$, $\tilde\theta$, $\tilde\lambda$ are unknown. Infer the values of $x(t_0), \tilde\theta(\theta), \tilde{\lambda}$ from the measurements of $y(t;t_0,x_0,\tilde{\theta},\tilde{\lambda})=C^{T}x(t;t_0,x_0,\tilde{\theta},\tilde\lambda)$ over $[t_0,t_0+T]$.
\end{prblm}

 The question is if there is an equivalent integral formulation such as e.g. (\ref{eq:integral_formulation_general}) of this problem for (\ref{eq:system1})? If such an integral formulation exists then whether a reduced-complexity version of this formulation can be stated so that the dimension of the parameter vector in the reduced formulation is smaller than that of in the original problem? Answers to these questions are provided in the next section.

\setcounter{equation}{0}

\section{Main Result}\label{sec:main_results}

\subsection{Indistinguishable parameterizations of (\ref{eq:system1})}

We begin with the following property of linear systems regarding input detectability (cf \cite{Tyukin_2012:arx})
\begin{lmm}\label{lem:observer_inferrence}
Consider
\begin{equation}\label{eq:system_io}
\begin{split}
& \begin{array}{ll}
\dot{x}&=A x + u(t)+d(t),\\
y&=C^{T}x, \ x(t_0)=x_0, \ x_0\in\Real^n,
\end{array}
\end{split}
\end{equation}
where
\[
A=\left(\begin{array}{cc} \begin{array}{c}
                                a_1\\
                                \vdots\\
                                a_n
                          \end{array} & \begin{array}{c}
                                          I_{n-1}\\
                                           \\
                                            0
                                            \end{array}
         \end{array}\right), \ C=(1,0,\dots,0)^{T},
\]
and
$x, u,d:\Real\rightarrow\Real^n$, $u\in\mathcal{C}^1$,
$d\in\mathcal{C}$. Let $u(\cdot),\dot{u}(\cdot),d(\cdot)$ be bounded:
$\max\{\|u(t)\|,\|\dot{u}(t)\|\}\leq B, \
\|{d}(t)\|\leq \Delta_{\xi}$ for all $t\geq t_0$. Then the following hold:

\begin{itemize}
\item[1) ] if the solution of (\ref{eq:system_io}) is globally bounded for all $t\geq t_0$ then, for $T$ sufficiently large, there are $\kappa_1,\kappa_2\in\mathcal{K}$:
\[
\begin{split}
&\|y(\tau)\|_{\infty,[t_0,t_0+T]}\leq \varepsilon \Rightarrow  \ \exists \ t'(\varepsilon,x_0)\geq t_0: \ \left\|z_1(\tau) + u_1(\tau) \right\|_{\infty,[t',t_0+T]}\leq \kappa_1(\varepsilon) + \kappa_2(\Delta_{\xi}),
\end{split}
\]
where $z_1=(1,0,\dots,0) z$,
\begin{equation}\label{eq:system_io_filter}
\begin{split}
\dot z&= \Lambda  z +   G  u,\ \Lambda=\left(\begin{array}{ccc}
- b \  & \begin{array}{c}\vdots\\ \vdots \end{array} &
\begin{array}{c} I_{n-2}\\ 0
\end{array}\end{array}\right),\\
 G&=\left(\begin{array}{cc} - { b} & I_{n-1} \end{array}\right), \  z(t_0)=0,
\end{split}
\end{equation}
and  $ b=(b_1,\dots,b_{n-1})^{T}$:  real parts of the roots
of $s^{n-1}  + b_1 s^{n-2}+\cdots + b_{n-1}$ are negative.

\item[2)] if $ d(t)\equiv 0$, then $y(t)= 0$ for all $t\in[t_0,t_0+T]$ implies existence of $p\in\Real^{n-1}$
\begin{equation}\label{eq:identifiability}
(1,0,\dots,0)e^{\Lambda(t-t_0)} p+z_1(t)+u_1(t)=0
\end{equation}
for all  $t\in[t_0,t_0+T]$.
\end{itemize}
\end{lmm}
Proof of Lemma \ref{lem:observer_inferrence} is provided in the Appendix.

According to Lemma \ref{lem:observer_inferrence} the following two sets of parameters, associated with every $\tilde\theta, \tilde\lambda$, need special consideration. The first set is defined as
\[
\begin{split}
&\mathcal{E}_0(\tilde{\theta},\tilde{\lambda},T)=\{(\theta',\lambda'), \ \theta'\in\Real^r,\lambda'\in\Real^{d+k} \ | \\
& b \varphi(y(t),t)^{T} (\theta'-\tilde\theta) + g(y(t),{\lambda}',t)-g(y(t),\tilde{\lambda},t)=0  \ \mbox{for all} \ t\in[t_0,t_0+T] \}.
\end{split}
\]
The set $\mathcal{E}_0(\tilde{\theta},\tilde{\lambda},T)$ contains all parameterizations of (\ref{eq:system1}) which are indistinguishable from each other providing that the values of  $x(t)$ are known for all $t\in[t_0,t_0+T]$. That is, if $x(t;t_0,x_0,\tilde\theta,\tilde\lambda)=x(t;t_0,x_0,\theta',\lambda')$ for all $t\in[t_0,t_0+T]$ then $(\theta',\lambda')\in\mathcal{E}_0(\tilde\theta,\tilde\lambda,T)$.

Denote
\[
\begin{split}
&\eta(\tilde\theta,\tilde\lambda,\theta',\lambda',p,t)=\varphi(y(t),t)^{T} (\theta'-\tilde\theta) +  g_1(y(t),{\lambda}',t)\\
&-g_1(y(t),\tilde{\lambda},t)+\tilde{C}^T e^{\Lambda(t-t_0)}p + z_1(t;t_0,\lambda')-z_1(t;t_0,\tilde\lambda),
\end{split}
\]
where $\Lambda,\tilde C, z(t;t_0,\lambda')$ are defined as in (\ref{eq:system_io_filter}) with $u(t)$ replaced by $g(y(t),\lambda',t)$. The second set is defined as
\[
\begin{split}
&\mathcal{E}(\tilde\theta,\tilde\lambda,T)=\{(\theta',\lambda'), \ \theta'\in\Real^r,\lambda'\in\Real^{d+k} \ | \\ \exists \ & p(\tilde\theta,\tilde\lambda,\theta',\lambda')\in \Real^{n-1}: \ \eta(\tilde\theta,\tilde\lambda,\theta',\lambda',p,t)=0 \ \mbox{for all} \ t\in[t_0,t_0+T] \}.
\end{split}
\]
In accordance with Lemma \ref{lem:observer_inferrence} the set $\mathcal{E}(\tilde\theta,\tilde\lambda,T)$ contains all parametrization of (\ref{eq:system1}) that are indistinguishable on the interval $[t_0,t_0+T]$ on the basis of accessing only the values of $y(x(t;t_0,x_0,\theta,\lambda))$. In other words, if $y(x(t;t_0,x_0,\tilde\theta,\tilde\lambda))=y(x(t;t_0,x_0',\theta',\lambda'))$ for all $t\in[t_0,t_0+T]$ then $(\theta',\lambda')\in\mathcal{E}(\tilde\theta,\tilde\lambda,T)$. If the set $\mathcal{E}(\tilde\theta,\tilde\lambda,T)$ contains more than one element then (\ref{eq:system1}) is not uniquely identifiable  on $[t_0,t_0+T]$ \cite{Distefano:1980}. Here, for simplicity, we will focus on systems (\ref{eq:system1}) that are uniquely identifiable on $[t_0,t_0+T]$:
\begin{assumption}\label{assume:unique_ident} Sets $\mathcal{E}_0(\tilde\theta,\tilde\lambda,T)$ and $\mathcal{E}(\tilde\theta,\tilde\lambda,T)$ coincide and contain no more than one element.
\end{assumption}

\subsection{Integral reduced-order formulation of the inverse problem for (\ref{eq:system1})}

Before we proceed with presenting an equivalent integral formulation of Problem \ref{prob:inferrence} let us first introduce several additional components and corresponding technical assumptions. Let $l\in\Real^n$ be a vector satisfying the following condition:
\[
P(A_0+l C^T)+(A_0+l C^T)^T P = -Q, \ Pb=C,
\]
where $P,Q$ are some symmetric positive definite matrices. According to the Meyer-Kalman-Yakubovich-Popov lemma, such vector will always exist since the polynomial $s^{n-1}+b_1 s^{n-2}+ \cdots + b_{n-1}$ is Hurwitz.

Consider
\begin{equation}\label{eq:error_dynamics}
\frac{d}{dt}\left(\begin{array}{c} \xi_1\\ \xi_2\end{array}\right)=\left(\begin{array}{cc} A_0+ l C^T & b\varphi(y(t),t)\\
                                    - \varphi(y(t),t) C^T & 0 \end{array}\right) \left(\begin{array}{c} \xi_1\\ \xi_2\end{array}\right),
\end{equation}
and let $\Phi(t,t_0)$ be its corresponding normalized fundamental solutions matrix: $\Phi(t_0,t_0)=I_{n+r}$.

\begin{thrm}\label{theorem:integral_formulation} Consider (\ref{eq:system1}) and suppose that Assumption \ref{assume:unique_ident} holds. Let $y(\cdot)$, $\varphi(y(\cdot),\cdot)$, $g(y(\cdot),\lambda,\cdot)$ be $T$-periodic on $[t_0,\infty]$ for all $\lambda$, and the function $\varphi(y(\cdot),\cdot)$ satisfy:
\[
\int_{t_0}^{t_0+T}\varphi(y(\tau),\tau) \varphi(y(\tau),\tau)^{T}d\tau \geq \delta I_r, \ \delta >0.
\]

Then the following statements are equivalent
\begin{itemize}
\item [1)]  $\hat{y}(\lambda',t)=y(t)$ for all $t\in [t_0,t_0+T]$, where $\hat{y}:\Real^{d+k}\times\Real\rightarrow\Real$:
\begin{equation}\label{eq:integral_formulation_canonic}
\begin{split}
&\hat{y}(\lambda',t)=(1 \ 0 \ \dots \ 0)\big(\Phi(t,t_0)R(\lambda')+\Phi(t,t_0)\int_{t_0}^{t}\Phi(\tau,t_0)^{-1} \left(\begin{array}{c} g(y(\tau),\lambda',\tau)-ly(\tau)\\ y(\tau)\varphi(y(\tau),\tau) \end{array}\right)d\tau \big)\\
&R(\lambda')=(I_{n+r}-\Phi(t_0+T,t_0))^{-1}\Phi(t_0+T,t_0)\int_{t_0}^{t_0+T}\Phi(\tau,t_0)^{-1} \left(\begin{array}{c} g(y(\tau),\lambda',\tau)-l y(\tau)\\ y(\tau)\varphi(y(\tau),\tau) \end{array}\right)d\tau.
\end{split}
\end{equation}
\item [2)] $(1 \ 0  \ \cdots \ 0)x(t;t_0,x_0,\tilde\theta,\lambda')=y(t)$ for all $t\in [t_0,t_0+T]$.
\end{itemize}

Furthermore, the values of $x_0$, $\tilde\theta$ satisfy
\begin{equation}\label{eq:initial_conditions}
\left(\begin{array}{c}x_0\\
                    \tilde\theta
      \end{array}
\right)=R(\lambda').
\end{equation}

\end{thrm}

\begin{proof} Let us first show that 1) $\Rightarrow$ 2). Recall (see e.g. \cite{Lorea_2002}) that assumptions of the theorem imply  existence of positive numbers $\rho,D>0$:
\[
\|\Phi(t,t_0')\|\leq D e^{-\rho(t-t_0')} \ \mbox{for all} \ t\geq t_0',  \ t,t_0'\in[t_0,\infty).
\]
Hence there are no zero eigenvalues of the matrix  $I_{n+r}-\Phi(t_0+T,t_0)$, and $(I_{n+r}-\Phi(t_0+T,t_0))^{-1}$ exists.

Consider $\chi=(\chi_1,\chi_2)$:
\begin{equation}\label{eq:estimates}
\begin{split}
& \frac{d}{dt}\left(\begin{array}{c} \chi_1\\ \chi_2\end{array}\right)= \left(\begin{array}{cc} A_0+ l C^T & b\varphi(y(t),t)\\
                                    - \varphi(y(t),t) C^T & 0 \end{array}\right) \left(\begin{array}{c} \chi_1  \\ \chi_2\end{array}\right)  + \left(\begin{array}{c} g(y(t),\lambda',t)-ly(t)\\ y(t)\varphi(y(t),t) \end{array}\right)
\end{split}
\end{equation}
It is clear that solutions of (\ref{eq:estimates}) are defined for all $t\geq t_0$ providing that the definition of $y(\cdot)$, $g(y(\cdot),\lambda',\cdot)$, and $\varphi(y(\cdot),\cdot)$ are extended (periodically) on the interval $[t_0,\infty)$. Introduce the function $\zeta(\cdot)=(x(\cdot,t_0,x_0,\tilde\theta,\tilde\lambda),\tilde\theta)$ (in which the domain of the function $x(\cdot,t_0,x_0,\tilde\theta,\tilde\lambda)$ definition is extended to $[t_0,\infty)$), and consider the difference
\[
\xi=\chi-\zeta.
\]
Dynamics of $\xi$ satisfy (\ref{eq:error_dynamics}) with $\xi_1(t_0)=\chi_1(t_0)-x(t_0)$, $\xi_2(t_0)=\chi_2(t_0)-\tilde\theta$. Moreover, $\hat{y}(\lambda',t)=C^{T}\chi_1(t)$ for all $t \in[t_0,t_0+T]$ (or in $[t_0,\infty)$ if $\hat{y}(\lambda',\cdot)$ is periodically extended on $[t_0,\infty)$).

Let $\hat{y}(\lambda',t)\equiv y(t)$. This implies that $\chi_2-\tilde\theta = \mbox{const}$ for all $t\in[t_0,t_0+T]$. Hence according to Lemma \ref{lem:observer_inferrence} $(\chi_2(t_0),\lambda')$ belong to $\mathcal{E}(\tilde\theta,\tilde\lambda,T)$. Given that sets $\mathcal{E}(\tilde\theta,\tilde\lambda,T)$ and $\mathcal{E}_0(\tilde\theta,\tilde\lambda,T)$ coincide and contain just one element, $\tilde\theta,\tilde\lambda$, we conclude that $\chi_2(t_0)=\tilde\theta$, $\lambda'=\tilde\lambda$.

Notice that  $\lim_{t\rightarrow \infty}\xi(t)=0$ for all $\chi(t_0)$, and that
\begin{equation}\label{eq:stable_periodic}
\Phi(t,t_0)R(\lambda')+\Phi(t,t_0)  \int_{t_0}^{t}\Phi(\tau,t_0)^{-1} \left(\begin{array}{c} g(y(\tau),\lambda',\tau)-ly(\tau)\\ y(\tau)\varphi(y(\tau),\tau) \end{array}\right)d\tau
\end{equation}
is the unique exponentially stable periodic solution of (\ref{eq:estimates}). This implies that (\ref{eq:initial_conditions}) holds.

Let us show that 2) $\Rightarrow$ 1). Let $\tilde\theta,\lambda'$ be parameters for which the following identity holds: $y(x(t;t_0,x_0,\tilde\theta,\lambda'))=y(t)$ for all $t\in[t_0,t_0+T]$. Consider the function $\zeta(\cdot)$ defined earlier.  Given that (\ref{eq:stable_periodic}) is the unique exponentially stable periodic solution of (\ref{eq:estimates}),  that $\lim_{t\rightarrow\infty}\zeta(t)=0$ for arbitrary choice of initial conditions (i.e. vectors $\tilde\theta$, $x(t_0)$, and $\chi_1(t_0)$, $\chi_2(t_0)$) and that $\zeta(t)\equiv 0$ if $\chi_1(t_0)=x_0$, $\chi_2(t_0)=\tilde\theta$,  one concludes that $\hat{y}(\lambda',t)=y(x(t;t_0,x_0,\tilde\theta,\lambda'))=y(t)$ for all $t\in[t_0,t_0+T]$.
\end{proof}

\begin{rmrk} One may argue that it is, in principle, possible to obtain integral formulations of the corresponding inverse problem without using adaptive observer-inspired structures. Note, however, that since the original matrix $A(\theta)$ is allowed to depend on unknown parameters $\theta$, explicit expressions of solutions of (\ref{eq:system}) will involve extra nonlinearly parameterized terms, $e^{A(\theta)(t-t_0)}$. If closed-form expressions are applied to (\ref{eq:system1}) then the drawback is that the overall unknown parameters vector is $(x_0,\tilde\theta,\tilde\lambda)$, and its dimension is $n+r+d+k$. In the proposed solution  dimension of the unknown parameters vector is reduced to $d+k$ which is advantageous for systems with large number of unknowns.
\end{rmrk}

\begin{rmrk} The uncertainty reduction achieved in the proposed method is due to the assumption that all functions in the right-hand side of (\ref{eq:system1}) are $T$-periodic. Whereas such periodicity assumptions may not always hold, they are not particularly difficult to satisfy (at least approximately) in the laboratory conditions.
\end{rmrk}

\begin{rmrk}\label{rem:discrete} Instead of dealing with continuous-time signals, $y(t)$, one may re-formulate the above results in the setting in which model responses and data are compared at mere $N$ discrete points $\{t_i\}$ in $[t_0,t_0+T]$.  In this case sets $\mathcal{E}_0$, $\mathcal{E}$ will need to be re-defined so that the corresponding identities hold at a finite number of points $\{t_i\}$ rather than for all $t\in[t_0,t_0+T]$. Discrete extension of the theorem allows straightforward formulation of the inference problem as
\begin{equation}\label{eq:discrete}
\tilde\lambda=\arg \min_{\lambda\in\Real^r} \sum_{i=1}^N (\hat{y}(\lambda,t_i)-y(t_i))^2
 \end{equation}
which bears some similarity with \cite{Kuhl_2011,Abarbanel:2009}. Here, however, no discretization of the original continuous-time dynamical model is required and  ${\pd \hat y(\lambda,t_i)}/{\pd \lambda}$ are computable as definite integrals.

Note also that in some applications measured data, $y(t)$, may only be available at certain discrete points. In this situation one can employ a suitable interpolation scheme providing that the outcome of such an interpolation is phenomenologically adequate.
\end{rmrk}

\begin{rmrk} Our method, as formulated, requires periodicity of $y(\cdot)$. Similar representations can be obtained for models that do not necessarily  produce periodic signals. This, however, will bear additional costs. In absence of periodicity  $R(\lambda')$ in (\ref{eq:integral_formulation_canonic}) will generally be replaced by an unknown vector. Yet, if the interval of observation is long enough then relative contribution of this unknown term in $y(\lambda',t)$ for $t$ sufficiently large will be small. Thus $y(\lambda',t)$ becomes an approximation of the measured $y(t)$ rather than an exact match.
\end{rmrk}

\begin{rmrk}
The method can also be linked with \cite{Pavlov_2013} where periodicity of input is used for inferring partial derivatives of outputs with respect to parameters and initial conditions. In \cite{Pavlov_2013} the main focus was on deriving a method for fast numerical estimation of unique periodic solutions of the $\delta_p$-subsystem in (\ref{eq:sensitivity}). Here we used a somewhat different approach whereby the problem is considered through the prism of adaptive observer design enabling us to express explicitly the values of measured trajectories as integrals of computable functions of model parameters.
\end{rmrk}

\setcounter{equation}{0}

\section{Examples}\label{sec:example}

\subsection{Predator-Prey system}

Consider (\ref{eq:predator_prey}), and let $x$ be the variable that is available for direct observation. Suppose that parameter values of (\ref{eq:predator_prey}) correspond to the unique stable limit cycle, and for simplicity we assume that the system evolves on the cycle (or in its sufficiently small vicinity). The corresponding parameter values and initial conditions are set as follows:
\begin{equation}\label{eq:predator_prey_parameters}
p_1=1, \ p_2=1.3, \ p_3=1, \ p_4=1, \ p_5=3, \ p_6=0.1, \ (x_0,z_0)=(0.0053,0.2536).
\end{equation}
Note that system (\ref{eq:predator_prey}) is not in the class of systems specified by (\ref{eq:system1}) to which Theorem \ref{theorem:integral_formulation} applies. It may, however, be transformed into this class as follows.
Given that $p_3, p_5\neq 0$ we will first introduce a new variable
\[
q=x + \frac{p_3}{p_5} z.
\]
Thus, in accordance with (\ref{eq:predator_prey}):
\[
\dot{q}=p_1 x \left(1 - \frac{x}{p_2} \right) - \frac{p_3 z x}{p_4 + x}+ \frac{p_3 z x}{p_4 + x} - p_6 \frac{p_3}{p_5}z,
\]
 and the system equations in the new coordinates become:
\begin{equation}\label{eq:predator_prey:modified}
\begin{split}
\dot{x}&= p_1 x \left(1 - \frac{x}{p_2} \right) - \frac{p_5\left(q-x\right) x}{\left(p_4  + x \right)}\\
\dot{q}& = p_1 x \left(1 - \frac{x}{p_2} \right) + p_6 x  - p_6 q.
\end{split}
\end{equation}
Equation (\ref{eq:predator_prey:modified}) is in the form (\ref{eq:1}). The latter, in turn, can be transformed into  (\ref{eq:system1}) by means of the following closed form expression for the variable $q$:
\[
q(t,p_1,p_2,p_6)= e^{-p_6(t-t_0)} q_0 (p_1,p_2,p_6)+ e^{-p_6 t}\int_{t_0}^t e^{p_6 \tau}\left[p_1x(\tau)\left(1-\frac{x(\tau)}{p_2} \right) + p_6 x(\tau)\right]d\tau.
\]
The observed variable, $x(\cdot)$, is periodic with period $T=34.05$, and $p_6=0.1 \neq 0$. Therefore
\[
q_0(p_1,p_2,p_6)=(1-e^{-p_6 T})^{-1} e^{- p_6 T} \int_{t_0}^{t_0+T} e^{p_6 \tau}\left[p_1x(\tau)\left(1-\frac{x(\tau)}{p_2} \right) + p_6 x(\tau)\right]d\tau,
\]
Hence dynamics of $x$ obeys
\begin{equation}\label{eq:example_predator_reduced}
\dot{x}= p_1 x - \frac{p_1}{p_2} x^2 + \frac{p_5 x^2}{\left(p_4  + x \right)} - \frac{p_5 x}{\left(p_4  + x \right)} q(t,p_1,p_2,p_6),
\end{equation}
which is of  class (\ref{eq:1}) or (\ref{eq:system1}) with
\[
\tilde{\lambda}=(p_1,p_2,p_4,p_5,p_6).
\]
Thus Theorem \ref{theorem:integral_formulation} applies, and observed periodic trajectory $x(\cdot)$ of system (\ref{eq:predator_prey}) can be represented as an explicit integral.

\begin{figure}
\includegraphics[width=0.45\linewidth]{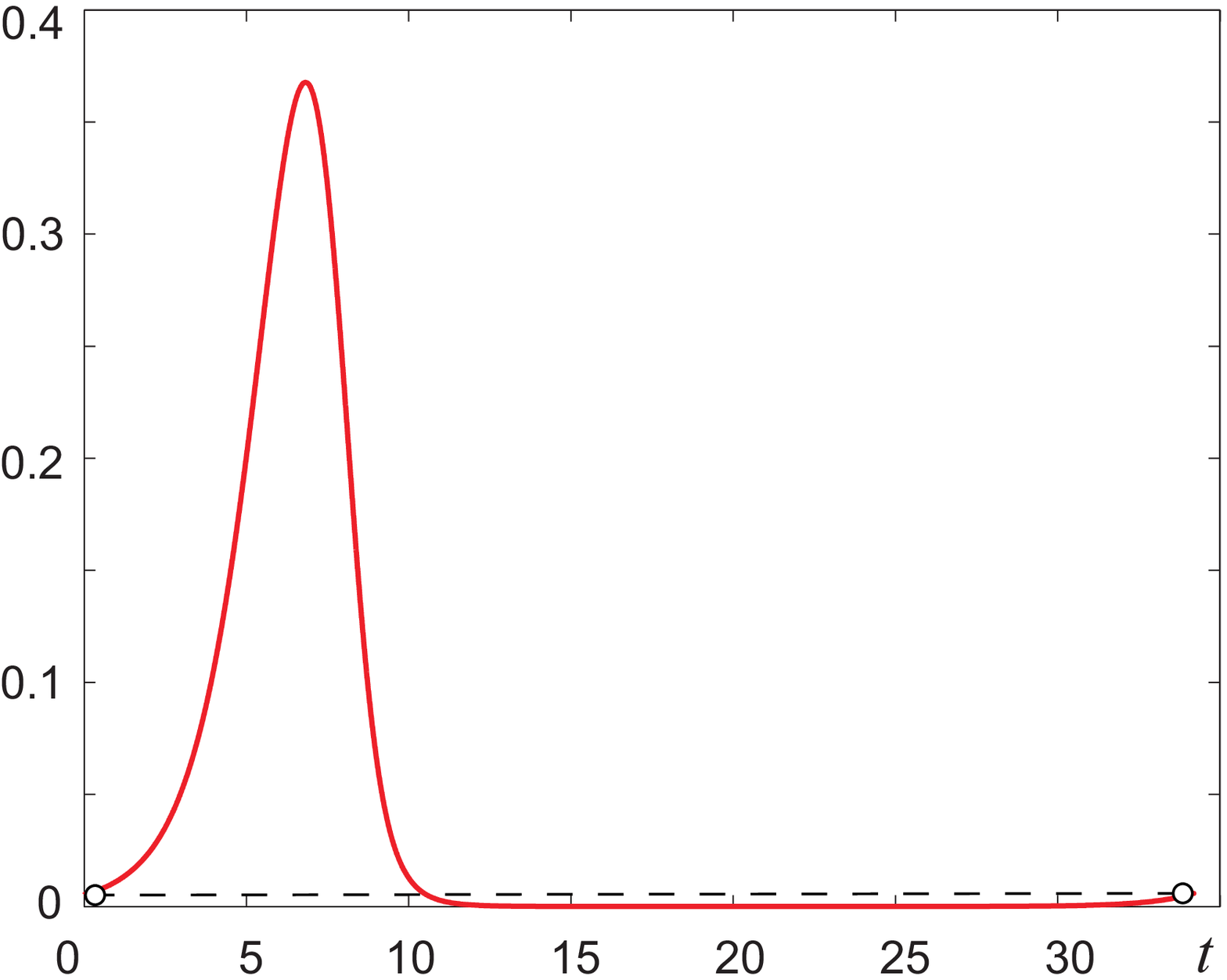}
\hspace{5mm}
\includegraphics[width=0.45\linewidth]{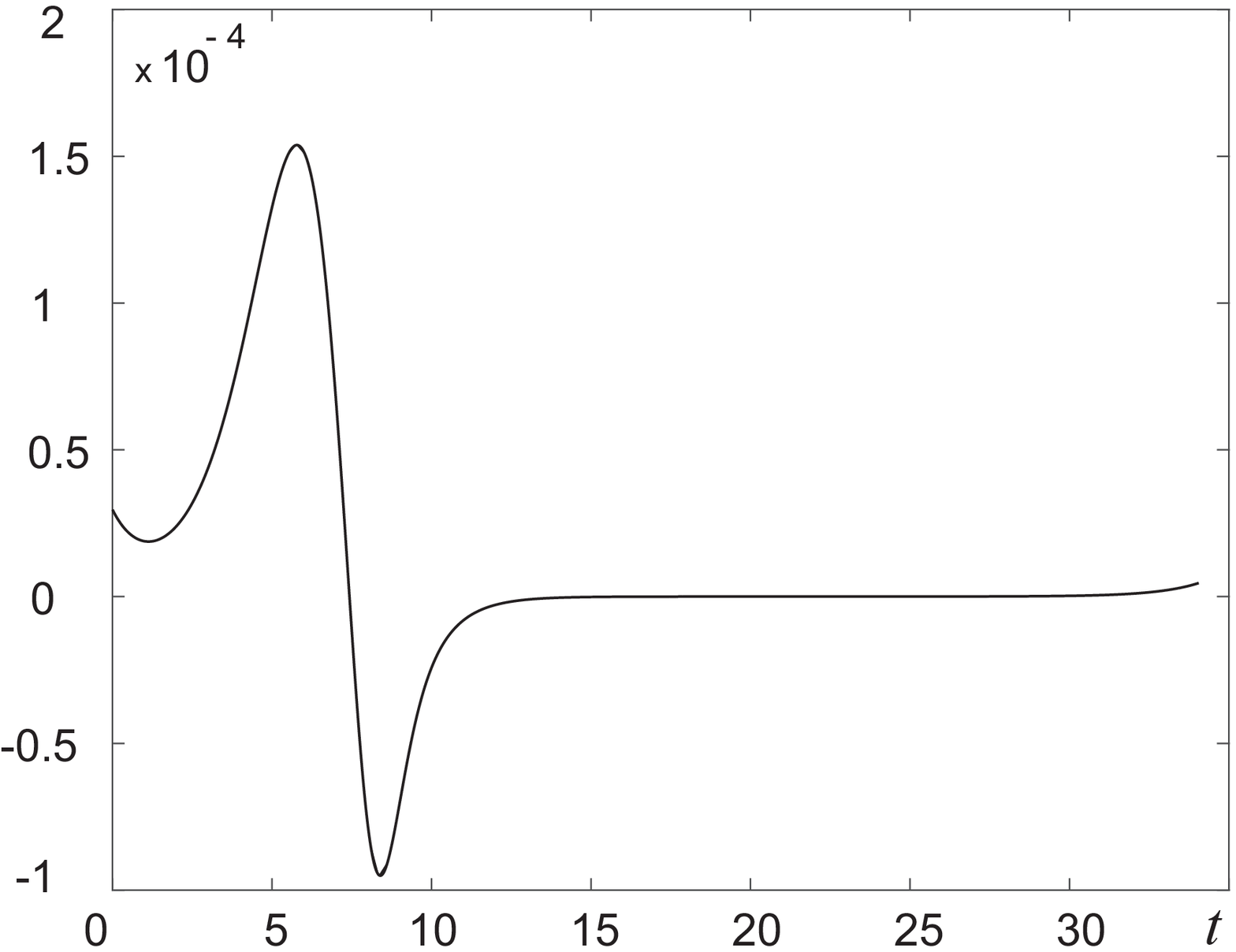}
\caption{Left panel: the values of $x(t;t_0,(x_0,z_0),p)$ and $\hat{y}(\tilde\lambda,t)$ as functions of $t$ for  $p=(p_1,...,p_6)$ and initial conditions specified by (\ref{eq:predator_prey_parameters}). Black circles indicate starting and ending points of the periodic trajectory $x(t;t_0,(x_0,z_0),p)$. The values of $x(t;t_0,(x_0,z_0),p)$ (red curve) were obtained by numerical integration of (\ref{eq:predator_prey}) by improved Euler with integration step $0.001$. The values of $\hat{y}(\tilde{\lambda},t)$ (blue curve) have been computed from representation (\ref{eq:predator_prey_integral}) numerically by simple right-hand rectangular integration with the same integration step. Right panel: the values of error, $\hat{y} (\tilde\lambda,t)-x(t;t_0,(x_0,z_0),p)$ as a function of $t$.}\label{fig:predator_prey}
\end{figure}

Notice that the number of parameters in (\ref{eq:example_predator_reduced}) is reduced to just $5$ as compared to $6$ in the original equations. Moreover, since the right-hand side of (\ref{eq:example_predator_reduced}) is purely nonlinearly parameterized, there is no $\varphi(\cdot)$ in (\ref{eq:error_dynamics}), $A_0=0$, and the fundamental solution matrix, $\Phi(t,t_0)$, becomes
\[
\Phi(t,t_0)=e^{l(t-t_0)}, \ l\in\Real, \ l<0.
\]
For the sake of simplicity we set $l=-1$. The corresponding expression for $\hat{y}(\tilde\lambda,t)$ is
\begin{equation}\label{eq:predator_prey_integral}
\begin{split}
\hat{y}(\tilde\lambda,t)=&e^{-(t-t_0)} R(\tilde\lambda) + \int_{t_0}^{t} e^{-(t-\tau)}\left(x(\tau) + p_1 x(\tau) - \frac{p_1}{p_2} x^2(\tau) + \frac{p_5 x^2(\tau)}{\left(p_4  + x(\tau) \right)}\right. \\
& \left.\ \ \ \ \ \ \ \ \ \ \ \ \ \ \ \ \ \ \ \ \ \ \ - \frac{p_5 x(\tau)}{\left(p_4  + x(\tau) \right)} q(\tau,p_1,p_2,p_6)\right)d\tau,\\
R(\tilde\lambda)=&(1-e^{-T})^{-1}\int_{t_0}^{t_0+T} e^{-(T-\tau)}\left(x(\tau) + p_1 x(\tau) - \frac{p_1}{p_2} x^2(\tau) + \frac{p_5 x^2(\tau)}{\left(p_4  + x(\tau) \right)}\right. \\
& \left.\ \ \ \ \ \ \ \ \ \ \ \ \ \ \ \  \ \ \ \ \ \ \ - \frac{p_5 x(\tau)}{\left(p_4  + x(\tau) \right)} q(\tau,p_1,p_2,p_6)\right)d\tau.
\end{split}
\end{equation}
Trajectories $x(\cdot;t_0,(x_0,z_0),p)$ and $\hat{y}(\tilde\lambda,\cdot)$ are shown in Fig. \ref{fig:predator_prey}.
Notice that trajectories $x(\cdot;t_0,(x_0,z_0),p)$ and $\hat{y}(\tilde\lambda,\cdot)$ nearly coincide with discrepancies of the order of $10^{-4}$ that are due to the differences in numerical integration.

In order to illustrate how our method works for this class of systems let us suppose that $y(t_i)=x(t_i;t_0,(x_0,z_0),p)$ be the measured data, and values of $p_1,\dots,p_6$ be unknown. For this particular simulation $t_i$ formed an equispaced grid in $[t_0,t_0+T]$ with $t_{i+1}-t_i=0.001$ for all $i=0,1,\dots,N-1$. The values $x(t_i;t_0,(x_0,z_0),p)$ were derived using the Runge-Kutta $4$th order method. As a measure of closeness between $y(\cdot)$ and $\hat{y}(\tilde\lambda,\cdot)$ we used the sum $\sum_{i=0}^N (y(t_i)-\hat{y}(\tilde\lambda,t_i))^2$, and as a parameter estimation routine we used the Nelder-Mead algorithm \cite{NelderMead}. The values of reflection, expansion, and contraction coefficients in the algorithm were set to $1,2$, and $0.5$, respectively.

\begin{figure}
\centering
\includegraphics[width=0.65\linewidth]{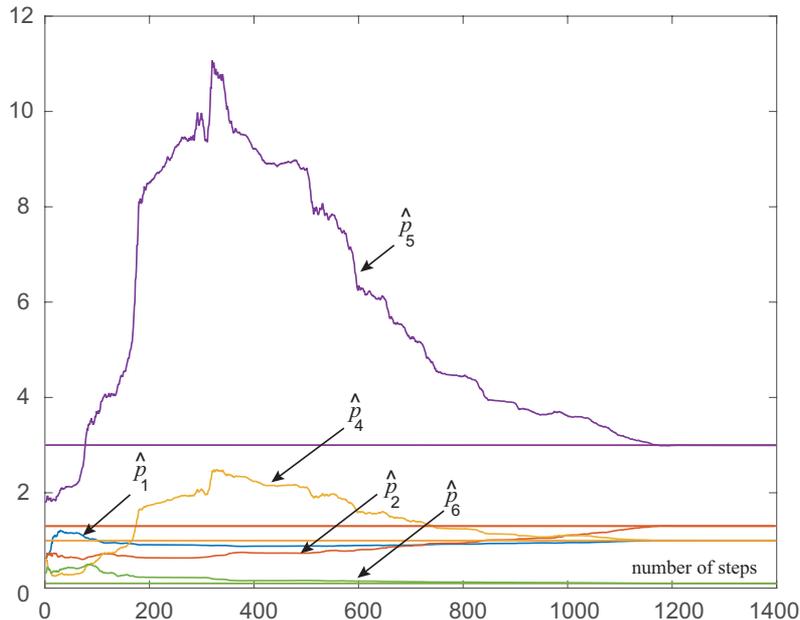}
\caption{Evolution of estimates $\hat{p}_1,\hat{p}_2,\hat{p}_4,\hat{p}_5,\hat{p}_6$ of true values of parameters $p_1,p_2,p_4,p_5,p_6$.}\label{fig:predator_prey_estimates}
\end{figure}

Behavior of parameter estimates are shown in Fig. \ref{fig:predator_prey_estimates}, and their initial and final values are provided in Table \ref{table:predator_prey_parameters}. As one can see from the table and plots, after roughly $1400$ steps the estimates are already in close proximity of true values of model parameters.
\begin{table}
\caption{True (first row), Initial (second row), and Estimated (third row) parameter values of (\ref{eq:predator_prey:modified}).}\label{table:predator_prey_parameters}
\begin{center}
Vector $\tilde\lambda=(p_1,p_2,p_4,p_5,p_6)$\\
\vspace{2mm}
\begin{tabular}{|c|c|c|c|c|}
\hline
 $p_1$ & $p_2$ & $p_4$ & $p_5$ & $p_6$\\
\hline
 $1$ & $1.3$ & $1$ & $3$ & $0.1$ \\
\hline
$0.3$ & $0.3$ & $0.3$ & $1.5$ & $0.01$\\
\hline
 $0.9999$ & $1.3018$ & $0.9991$ & $2.9966$ & $0.1$\\
\hline
\end{tabular}\\
\end{center}
\end{table}
This particular simulation took $21.06$ seconds in Matlab R2015a on standard desktop with Intel i5-2400 CPU and 8Gb of memory.

Let us now consider another relevant model that will enable us not only to demonstrate existence of explicit integral representation of its solutions but also to illustrate the point that sometimes the overall number of unknown parameters can be reduced as a result of the proposed integral representation.

\subsection{Morris-Lecar system}

Consider a model describing dynamics of voltage oscillations generated in barnacle giant muscle fiber \cite{Morris_Lecar}. The model can be viewed as a physiologically relevant reduction (see e.g. \cite{Izhikevich:2007}) of the celebrated Hodgkin-Huxley system \cite{Hodgkin_Huxley} of equations. The model equations are as follows:
\begin{equation}\label{eq:Morris-Lecar}
\begin{split}
\dot{x}=&g_{Ca}m_{\infty}(x)(x+E_{Ca})+g_{K}q(x+E_K)+g_L(x+E_L)+I\\
\dot{q}=&-\frac{1}{\tau(x)}q+\frac{w_\infty(x)}{\tau(x)},\\
y=&x,
\end{split}
\end{equation}
where
\[
\begin{split}
m_\infty(x)&=0.5\left(1+\tanh\left(\frac{x-V_1}{V_2}\right)\right)\\
w_\infty(x)&=0.5\left(1+\tanh\left(\frac{x+V_3}{V_4}\right)\right)\\
\tau(x)&=T_0 \left(\cosh\left(\frac{x+V_3}{2V_4}\right)\right)^{-1}
\end{split}.
\]
 Variable $x$ is the measured voltage, $q$ is the recovery variable. The values of $E_{Ca}$, $E_K$, $E_L$ are normally known ($E_{Ca}=-100$, $E_K=70$, $E_L=50$); other parameters may vary from one cell to another.

It is clear that equations (\ref{eq:Morris-Lecar}) are of the form (\ref{eq:system}). Moreover, if the model operates in the oscillatory regime then the right-hand side is periodic in $t$, including the variable $q$. In addition the integral
\[
\int_{t_0}^{t_0+T}-\frac{1}{\tau(x(s))}ds < 0,
\]
where if $T$ is the period of oscillations, for practically relevant values of $T_0,V_3,V_4$. Assuming that observations are taking place when the system's solution are on (or sufficiently near) the stable period orbit we can express the variable $q(t)$  as follows:
\[
\begin{split}
q(t)=&e^{\int_{t_0}^{t}-\frac{1}{\tau(x(s))}ds} q_0 + \int_{t_0}^{t} e^{\int_{z}^{t}-\frac{1}{\tau(x(s))}ds}\frac{w_\infty(x(z))}{\tau(x(z))}dz\\
q_0=& (1-e^{\int_{t_0}^{t_0+T}-\frac{1}{\tau(x(s))}ds})^{-1}\int_{t_0}^{t_0+T} e^{\int_{z}^{t_0+T}-\frac{1}{\tau(x(s))}ds}\frac{w_\infty(x(z))}{\tau(x(z))}dz.
\end{split}
\]
This brings equations (\ref{eq:Morris-Lecar}) into the form (\ref{eq:system1}) with parameters $\tilde\theta=(g_L,I)$, and $\tilde\lambda=(V_1,V_2,V_3,V_4,T_0,g_{Ca},g_K)$.

For the purpose of illustration we set the values of parameters $\tilde\theta$, $\tilde\lambda$ as specified in Table \ref{table:parameters}.
\begin{table}
\caption{True (first row) and Estimated (second) parameter values of (\ref{eq:Morris-Lecar}).}\label{table:parameters}
\label{table_example}
\begin{center}
Vector $\tilde\lambda=(V_1,V_2,V_3,V_4,T_0,g_{Ca},g_K)$\\
\vspace{2mm}
\begin{tabular}{|c|c|c|c|c|c|c|}
\hline
 $V_1$ & $V_2$ & $V_3$ & $V_4$ & $T_0$ & $g_{Ca}$ & $g_K$\\
\hline
 $-1$ & $15$ & $-10$ & $14.5$ & $3$ & $-1.1$ & $-2$\\
\hline
 $-0.95$ & $15.08$ & $-10.15$ & $14.44$ & $3.04$ & $-1.12$ & $-2.02$\\
\hline
\end{tabular}\\
\vspace{5mm}
Vector $\tilde\theta=(g_L,I)$\\
\vspace{2mm}
\begin{tabular}{|c|c|}
\hline
$g_L$  & $I$ \\
\hline
 $-0.5$  & $10$  \\
\hline
 $-0.539$  & $10.65$ \\
\hline
\end{tabular}
\end{center}
\end{table}
For the data generated at these parameter values the system is uniquely identifiable, and hence Assumption \ref{assume:unique_ident} holds. According to Theorem \ref{theorem:integral_formulation}, the problem of finding the values of $\tilde\theta,\tilde\lambda$ can be now formulated as that of matching the function $\hat{y}(\lambda',t)$ defined in (\ref{eq:integral_formulation_canonic}) to $y(t)$ over $[t_0,t_0+T]$. And in view of Remark \ref{rem:discrete} it reduces to
solving the unconstrained program (\ref{eq:discrete}).

In order to evaluate $\hat{y}(\lambda',t)$, as a function of parameter $\lambda'$ at a given $t$ one needs to know the fundamental solutions matrix $\Phi(t,t_0)$ for all $t\in[t_0,t_0+T]$. In this example this matrix was constructed numerically (using Dormand-Prince method and with fixed step size $0.0002$) from linearly independent solutions of
\begin{equation}\label{eq:example_fundamental}
\dot{z}=\left(\begin{array}{ccc} l  & y(t) & 1\\
                                   -y(t) & 0 & 0\\
                                   -1 & 0 & 0  \end{array}\right)z, \ l=-1
\end{equation}
starting from $(1,0,0)^{T}$, $(0,1,0)^{T}$, and $(0,0,1)^{T}$.

\begin{figure}
\centering
\includegraphics[width=0.5\linewidth]{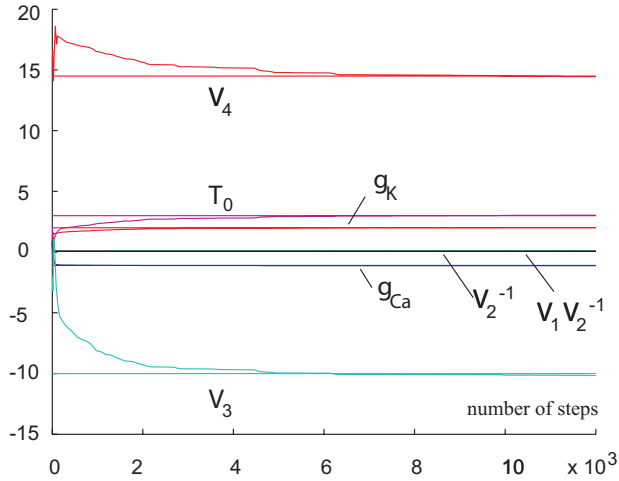}
\caption{Estimates and true values of $V_1/V_2$, $1/V_2$, $V_3$, $V_4$, $T_0$, $g_{Ca}$, $g_K$.}\label{fig: estimates}
\end{figure}

\begin{figure}
\centering
\includegraphics[width=0.4\linewidth]{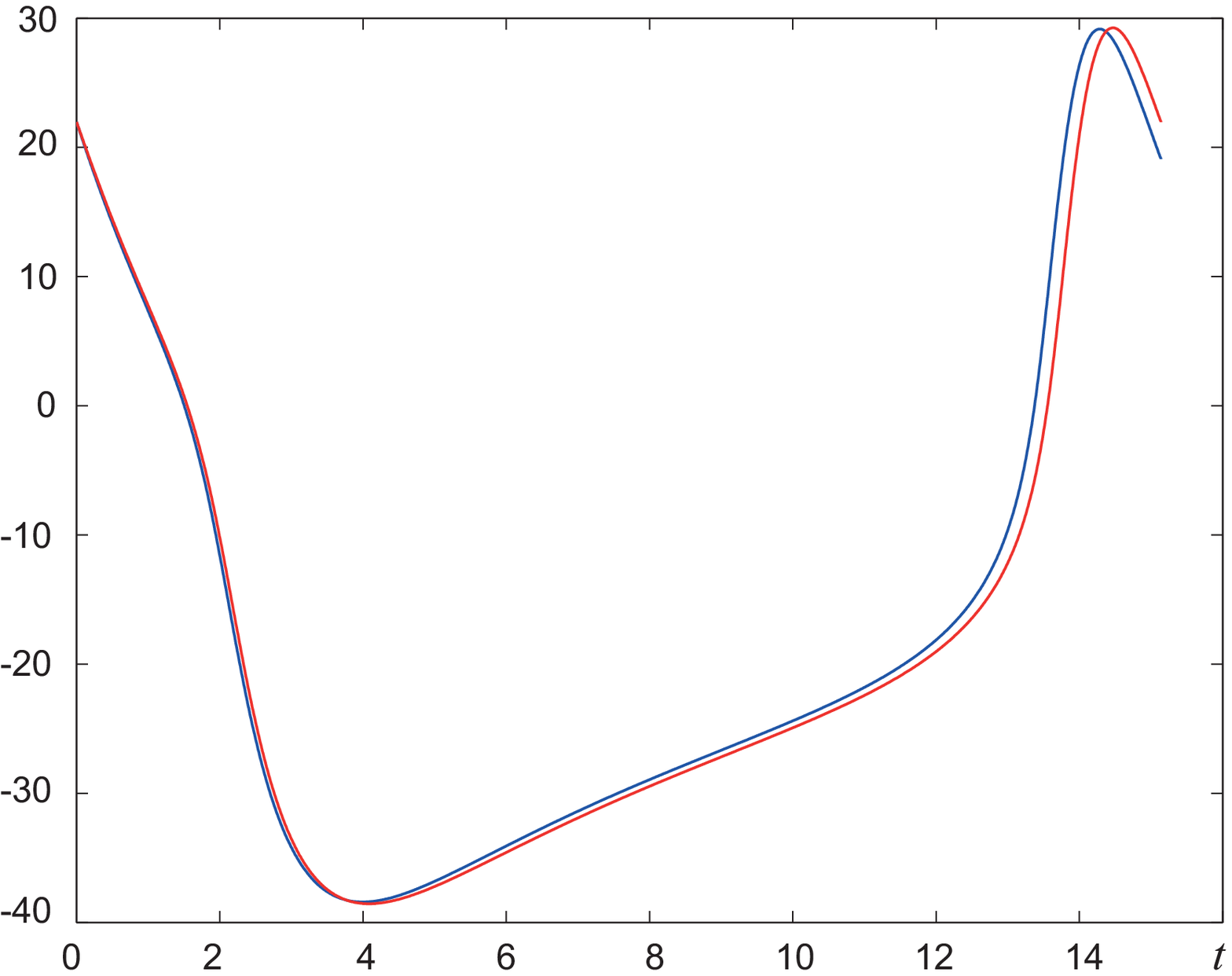}
\hspace{5mm}
\includegraphics[width=0.4\linewidth]{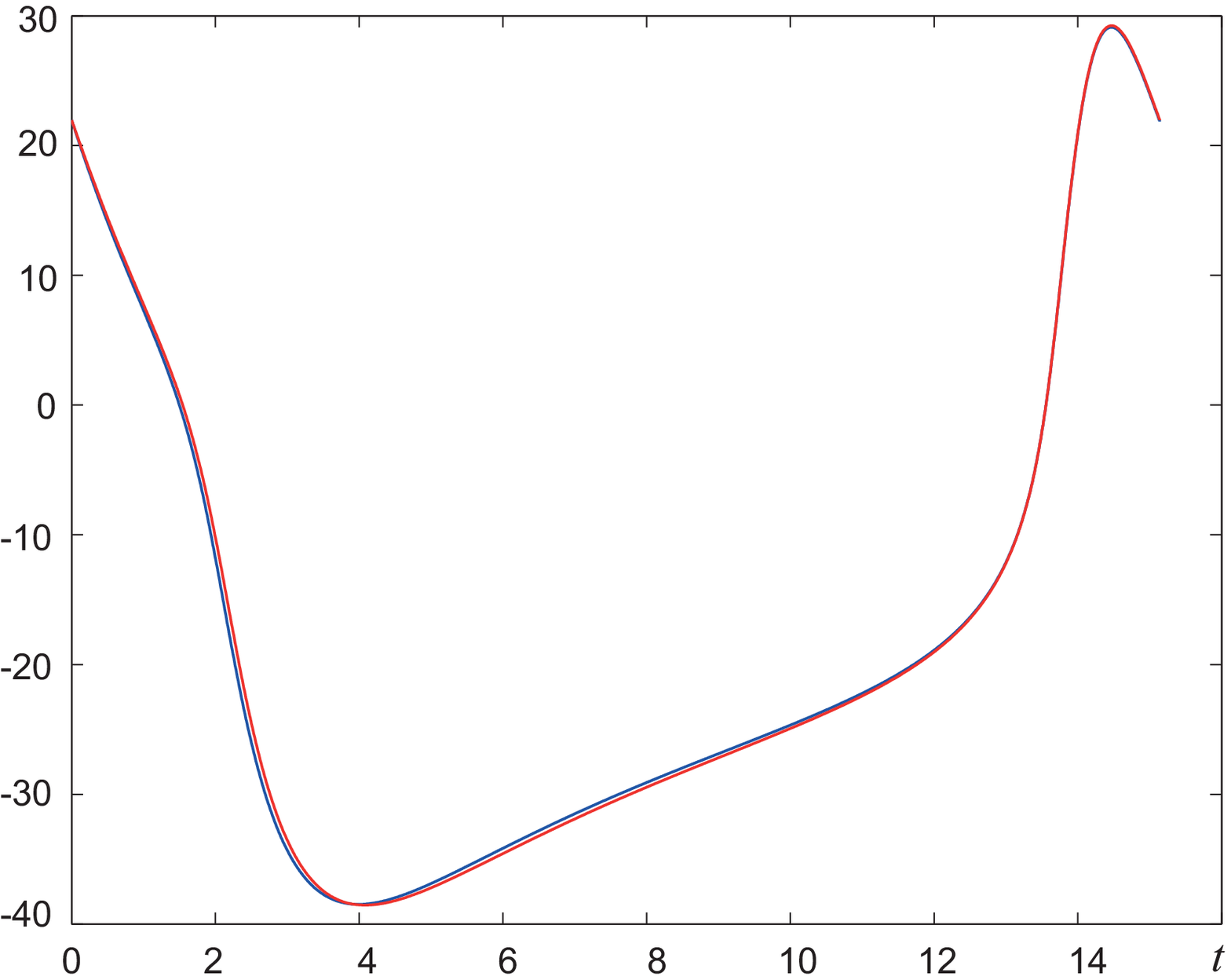}
\caption{Variables $x(t)$ of the state vector of (\ref{eq:Morris-Lecar}) for true values of $\tilde\theta,\tilde\lambda$ (red curves) and estimated values of $\tilde\theta,\tilde\lambda$ from Table \ref{table:parameters} (blue curves). The left panel shows the case when no adjustments to estimated parameters were made. The right panel illustrates how the reconstructed $x(t)$ changes when the parameter $\tilde{\theta}_2$ regulating the frequency of oscillations is slightly adjusted by $-0.07$.}\label{fig:M-L-trajectoris}
\end{figure}

\begin{figure}
\centering
\includegraphics[width=0.6\linewidth]{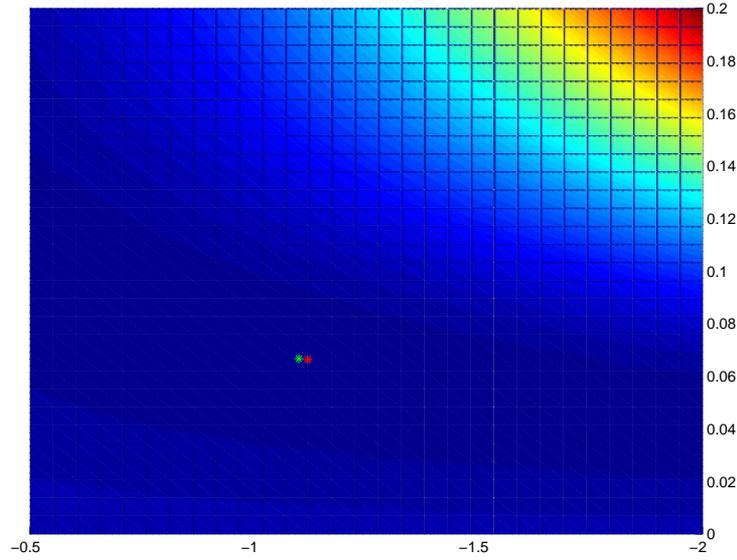}
\caption{Estimation error $\sum_{i=1}^N (\hat{y}(\tilde\lambda,t_i)-y(t_i))^2$ plotted as a function of $1/V_2$, $g_{Ca}$. Red star marks estimated $(g_{Ca},1/V_2)$, and green star corresponds to the true values of $(g_{Ca},1/V_2)$.  }\label{fig:M-L-ill}
\end{figure}

Points $t_i$ in (\ref{eq:discrete}) were evenly spaced with $t_{i+1}-t_i=0.04$, and the BFGS quasi-Newton method was used to find a numerical estimation of the solution of (\ref{eq:discrete}). For computational convenience, instead of looking for $V_1,V_2$ directly we were estimating ratios $1/V_2$ and $V_1/V_2$ respectively. Similarly, as follows from (\ref{eq:example_fundamental}), the estimate of parameter $I$ is not the value of $\tilde{\theta}_2$ but rather is the sum $\tilde{\theta}_2-\tilde{\theta}_1 50$.
We run the method for $12000$ iterations, and results  of the estimation are shown in Table \ref{table:parameters} and Fig. \ref{fig: estimates}.
In order to verify the quality of parameter estimation we run (\ref{eq:Morris-Lecar}) with both estimated and true values of parameters. Results of this simulation are show in Fig. \ref{fig:M-L-trajectoris}, upper panel. Note that frequency of the estimated $x(t)$ is higher than that of the measured data. This explains noticeable difference between trajectories at the end of the interval. In order to compensate for this difference we adjusted parameter $\tilde{\theta}_2$ (regulating the frequency of oscillations in the original model) by $-0.07$. Simulated trajectory of (\ref{eq:Morris-Lecar}) after this adjustment is shown in Fig. \ref{fig:M-L-trajectoris}, lower panel.
It is worth noticing that even though both estimated and simulated $x(t)$ are matching reasonably well there are still errors. The origin of these errors is likely to be 1) due to numerical errors in estimating the matrix $\Phi(t,t_0)$, and 2) due to the ill-conditioning of the original problem. Indeed, as Fig. \ref{fig:M-L-ill} suggests, there is a long shallow valley in a vicinity of the optimum.

The estimation took approximately $1$ hour on a standard PC (same configuration as in the previous example) in Matlab R2015a. We observed that most of the time was spent in the calculations of $\frac{\pd \hat{y}}{\pd \tilde\lambda}$ which is not surprising given the integration (\ref{eq:integral_formulation_canonic}) was performed over a relatively dense and uniform grid of points. On the other hand, this indicates that in this and similar cases scalability of the procedure is expected to grow nearly linear with dimension of $\tilde\lambda$. This will be tested in experiments in future.

In order to assess potential computational advantage of the proposed integral representations we compared the time needed for $1000$ evaluations of $y(t)$ on CPU and GPU over the interval $[t_0,t_0+T]$ for $1000$ randomly chosen parameter values. The results are summarized in Table \ref{tab:CPU_GPU_compare} (upper table). We see that GPU implementation of the same procedure resulted in the $39$-fold performance gain.
 Our second set of experiments assessed the time needed for running $1000$ explicit Euler integrations over the same period and for the same parameter values. The results are shown in Table \ref{tab:CPU_GPU_compare} (bottom table). These experiments showed that explicit integral representations, if implemented on GPU, are approximately $15$ times faster than direct Euler integration on CPU.

  We also notice that $1000$ model evaluations using explicit Euler integration on GPU is approximately $2$ times faster than the proposed integral implementations in this problem. Yet, our proposed scheme returns the estimates of all initial conditions and parameters that enter the right-hand side linearly ($2$ parameters, one initial condition). Furthermore, it enables to consolidate all computational power of the GPU into a single stream of computations which will be advantageous for local and inherently iterative optimization methods such as e.g. gradient-based search. In this regards comparing amount of time spent in integrating the corresponding sensitivity functions system with our explicit integral representation would be a fairer setting. This will be done in our future work.


\begin{table}
\begin{center}
\vspace{5mm}
Integral representation
\vspace{5mm}

\begin{tabular}{|c|c|c|}
\hline
CPU Xeon, $1000$ times  & GPU K10, $1000$ times & Ratio\\
\hline
$24.041$ sec  & $0.616$ sec & $39.009$ \\
\hline
\end{tabular}

\vspace{5mm}

Direct explicit Euler integration

\vspace{5mm}

\begin{tabular}{|c|c|c|}
\hline
CPU Xeon, Euler $1000$ times  & GPU K10, Euler $1000$ times & Ratio\\
\hline
$9.354$ sec  & $0.395$ sec & $23.65$ \\
\hline
\end{tabular}
\end{center}
\caption{Performance comparison chart.}\label{tab:CPU_GPU_compare}
\end{table}

\section{Conclusion}

We presented a technique for fast reconstruction of state and parameters from observed trajectories of evolving processes. The technique employs ideas from adaptive observers design to express measured trajectories as explicit functions of unknown parameters  and initial conditions. Such integral representation enables to exploit advantages of using parallel computational streams for efficient model evaluations, and presents a solution to a class of inverse problems that may be made scalable with computational resources available. The method applies to both linearly and nonlinearly parameterized systems. Moreover, it has been shown that the method allows to reduce dimensionality of the problem to that of the dimension of the vector of parameters entering the right-hand side of the model nonlinearly.

With respect to applicability of the method for predicting drastic changes in behavior of systems that are being investigated, we note that perhaps tracking only those bifurcations that persist under small perturbations is relevant for these purposes \cite{Gorban:1980}, \cite{Gorban:2004}. This approach to predictive modelling is similar to the ideas of structural sensitivity/insensitivity considered in  \cite{Adamson}.

The viability of the method was tested on two examples. The method performed well in these problems. We have also shown that the proposed approach can benefit from  parallel implementation of explicit numerical integration involved. It would be interesting to see if the same approach could be applied to a broader range of problems. Answering to this call will be the subject of our future work in this direction.

\begin{acknowledgement}
Ivan Tyukin is thankful to the Russian Foundation for Basic Research (research project No. 15-38-20178) for partial support.
\end{acknowledgement}

%
%

\setcounter{equation}{0}

\section{Appendix}\label{app:proofs}

\begin{lmm}\label{lem:system_io_inferrence_1order} Consider $\dot{y}=k y + u(t)+d(t)$, $k\in\Real$, $u,d:\Real_{\geq t_0}\rightarrow\Real$, $u\in\mathcal{C}^1$, $ d\in\mathcal{C}^{0}$, and let $\max\{|u(t)|,|\dot{u}(t)|\}\leq B$, $|d(t)|\leq \Delta_\xi$. Finally, let $T,\varepsilon$ be non-negative real numbers such that $T>\sqrt{\varepsilon}$. Then
\[
\begin{split}
&\|y\|_{\infty,[t_0,t_0+T]}\leq \varepsilon \Rightarrow   \|u\|_{\infty,[t_0,t_0+T)} \leq
\sqrt{\varepsilon}(1+e^{|k|\sqrt{\varepsilon}}+B)+\Delta_\xi.
\end{split}
\]
\end{lmm}

\begin{proof}
Let $L$ be an arbitrary element of $[0,T]$. Noticing that $y(t)$ for $t\geq t_0 + L$, $L>0$, can be expressed as:
$y(t)=y(t-L)e^{ k L}+\int_{t-L}^t
e^{k(t-\tau)}(u(\tau)+d(\tau))d\tau$ and using the Mean-value
theorem we obtain: $y(t)-y(t-L)e^{ k L}= L
e^{k(t-\tau')}(u(\tau')+d(\tau')), \ \tau'\in[t-L,t]$. Hence  $ \varepsilon(1+e^{kL})\geq L e^{k(t-\tau')}
(|u(t)|-L B-\Delta_{\xi})$, and
\[
\begin{split}
&\begin{array}{l}\Delta_\xi+LB+\frac{\varepsilon(1+e^{kL})}{L \min\{1,e^{kL}\}}\end{array} \geq \begin{array}{l}\Delta_\xi+LB+\frac{\varepsilon(1+e^{kL})}{L \min\{1,e^{k(t-\tau')}\}}\end{array}
\geq |u(t)|  \ \forall t\geq t_0+L.
\end{split}
\]
Given that $L$ can be chosen arbitrarily in the interval $[0,T]$ we let
$L=\sqrt{\varepsilon}$, and thus $|u(t)|\leq
\sqrt{\varepsilon}(1+e^{k\sqrt{\varepsilon}})\max\{1,e^{-k\sqrt{\varepsilon}}\}+B\sqrt{\varepsilon}
+\Delta_\xi \leq
\sqrt{\varepsilon}(1+e^{|k|\sqrt{\varepsilon}}+B)+\Delta_{\xi}
\ \forall \ t\in [t_0 + \sqrt{\varepsilon},t_0+T]$.

Finally, given that $|\dot{u}(t)|\leq B$ for all $t\in[t_0,t_0+T]$, including in the interval $[t_0,t_0+\sqrt{\varepsilon}]$, we conclude that
\[
|u(t)|\leq
\sqrt{\varepsilon}(1+e^{|k|\sqrt{\varepsilon}}+2B)+\Delta_{\xi}
\ \forall \ t\in [t_0,t_0+T].
\]

\end{proof}

\subsection{Proof of Lemma \ref{lem:observer_inferrence}}

Let us rewrite (\ref{eq:system_io}) as
\[
\begin{split}
\dot{y}&= a_{1} y + \tilde{C}\tilde x + u_1(t)+d_1(t)\\
\dot{\tilde x}&= \tilde A \tilde  x + \tilde  a y + b u_1 +  G  u(t)+\tilde{ d}(t),
\end{split}
\]
where $\tilde{a}=\mathrm{col}(a_2,\dots,a_n)$,
$\tilde C=\mathrm{col}(1,0,\dots,0)$,
$\tilde d(t)=\mathrm{col}(d_2(t),\dots,d_n(t))$, and
\[
G=\left(\begin{array}{cc} - {b} & I_{n-1}
\end{array}\right),  \
\tilde  A=\left(\begin{array}{cc}0 & I_{n-2}\\
                                0 & 0
                                \end{array}\right).
\]
Let $\|y(t)\|_{\infty,[t_0,t_0+T]}\leq \varepsilon$ and denote
$e(t)=\tilde C^{T}\tilde x+u_1(t)$.

According to Lemma
\ref{lem:system_io_inferrence_1order},  there are
$\upsilon_1,\upsilon_2\in\mathcal{K}$ such that
$\|e(t)\|=\|\tilde C^{T}\tilde x+u_1(t)\|\leq
\upsilon_1(\varepsilon)+\upsilon_2(\Delta_\xi)$ for all $t\in[t_0,t_0+T]$.

Using the notation above we obtain:
$\dot{\tilde x}=(\tilde  A -   b
\tilde C^{T})\tilde x + \tilde  a y(t) +
\tilde{ G} u(t) +   b e(t)+\tilde{ d}(t)$.

Matrix $\tilde  A -   b \tilde C^{T}=\Lambda$ is Hurwitz,
and hence there are $D,k\in\Real_{>0}$ such that
$\|e^{\Lambda(t-t_0)}\|\leq D e^{-k (t-t_0)}$.
Therefore $ \|\tilde  C^{T} \tilde x(t)-
\tilde C^{T}\int_{t_0}^{t}e^{\Lambda(t-\tau)} G
 u(\tau)d\tau\| \leq De^{-k(t-t_0)}\|\tilde{ x}(t_0)\|
+\frac{D}{k}(\| a\|\varepsilon
+\| b\|(\upsilon_1(\varepsilon)+\upsilon_2(\Delta_\xi))+\Delta_\xi)
$.

Noticing that
$z_1=\tilde C^{T}\int_{t_0}^{t}e^{\Lambda(t-\tau)} G
 u(\tau)d\tau$, denoting $\kappa(\varepsilon)=2
\frac{D}{k}(\| a\|\varepsilon
+\| b\|\upsilon_1(\varepsilon))+\upsilon_1(\varepsilon)$,
$\kappa_2(\Delta_\xi)=2
\frac{D}{k}(\Delta_\xi+\| b\|\upsilon_2(\Delta_\xi))+\upsilon_2(\Delta_\xi)$,
and
\[
t'(\varepsilon,x_0)= t_0 + \frac{1}{k}\ln\left(\frac{D\|x_0\|}{\varepsilon}\right)
\]
we can conclude that there is a $t'(\varepsilon,x_0) \geq t_0$ such that
\[
\begin{split}
\|z_1(\tau)+ u_1(\tau)\|_{\infty,[t,t_0+T]}&\leq \kappa(\varepsilon)+\varepsilon+\kappa_2(\Delta_\xi)=\kappa_1(\varepsilon)+\kappa_2(\Delta_\xi).
\end{split}
\]
for all $t\in [t'(\varepsilon,x_0),t_0+T]$, providing that $T$ is sufficiently large to satisfy $t_0+T > t'(\varepsilon,x_0)$.

Noticing that $y(t)\equiv 0 \Rightarrow e(t)\equiv0$ ensures that (\ref{eq:identifiability}) holds too. $\square$

\end{document}